\documentclass[12pt]{article}
\usepackage[english]{babel}
\usepackage{amssymb,amsmath,graphics,tikz}

\usepackage{hyperref}

\parskip\medskipamount
\newtheorem{theorem}{Theorem}[section]

\newtheorem{remark}[theorem]{Remark}

\newtheorem{lemma}[theorem]{Lemma}

\def\<{\langle}
\def\>{\rangle}
\newcommand{\proof}{\emph{Proof.~}}

\newcommand{\dd}{\mathsf{d}}

\newcommand{\PG}{\mathsf{PG}}

\newcommand{\cF}{\mathcal{F}}

\def\qed{{\hfill\hphantom{.}\nobreak\hfill$\Box$}}

\newcommand{\R}{\mathbb{R}}

\title{Descent of affine buildings - II. Minimal angle $\pi/3$ and exceptional quadrangles}
\author{ Koen Struyve\thanks{The author is supported by  the Fund for Scientific Research --
Flanders (FWO - Vlaanderen)} }
\begin{document}
\maketitle
\begin{abstract}
In this two-part paper we prove an existence result for affine buildings arising from exceptional algebraic reductive groups. Combined with earlier results on classical groups, this gives a complete and positive answer to the conjecture concerning the existence of affine buildings arising from such groups defined over a (skew) field with a complete valuation, as proposed by Jacques Tits.

This second part builds upon the results of the first part and deals with the remaining cases.
\end{abstract}
\section{Introduction}

As in the first part, we deal with the central problem to show the existence of the affine building associated with a semi-simple algebraic group
over a field with a complete valuation. In Part I we answered this question using descent methods positively when the `minimal angle' is strictly greater than $\pi/3$. Here, we will extend this result to the case where $B$ is of rank one and the minimal angle is exactly $\pi/3$.

While this might seem as only a small improvement,  it makes a huge difference for our main goal. If one computes the minimal angle of the Tits diagrams associated to the exceptional forms, they turn out to be almost all strictly greater than $\pi/3$ or they are of relative rank one and minimal angle equal to $\pi/3$, except for four cases. So for all but these four cases our geometric methods confirm Tits' conjecture. By adding some algebraic arguments we are also able to solve the remaining cases, including the one associated to the Moufang quadrangle of exceptional type $\mathsf{E}_8$.

We build upon the definitions and notations of Part I.



\section{Statement of the main result}

We now state the main result of the paper.


\textbf{Main Result.} \emph{
Let $(\Lambda,\cF)$ be an irreducible $\R$-building with a maximal system of apartments, and let $G$ be a bounded group of isometries acting on $\Lambda$. Suppose that the following holds:
\begin{itemize}
\item The action of $G$ on $\Lambda_\infty$ admits a Tits diagram of relative rank one whose minimal angle equals $\pi/3$.
\item $(\Lambda,\cF)$ is not a Bruhat-Tits building arising from a classical non-algebraic group.
\end{itemize}
Then one can construct an $\R$-tree $(\widetilde\Lambda,\widetilde\cF)$ such that
$\widetilde\Lambda_\infty$ is in a bijective correspondence with the fixed point set of $G$ in $\Lambda_\infty$.
}

%

This has the following corollary for the existence problem for affine buildings arising from exceptional algebraic groups.

\textbf{Main Corollary.}  \emph{
The conjecture proposed by Jacques Tits in~\cite[p. 173]{Tit:86} has a positive answer.
}

Details and proof for this corollary are given in Section~\ref{section:existence}.

\section{Proof of the main result}
The strategy of the proof is to start from the structure $\Lambda'$ obtained in Part I of this two-part paper and then modify it in a $\R$-tree via quotients.

\subsection{Additional definitions and notations}
We assume that we are in the situation as described by the statement of the main result. Let $(W,S)$ be the type of the $\R$-building $(\Lambda,\cF)$. 

The Tits diagram $(M,\Gamma,A)$ associated to the action of $G$ is of relative rank one, or equivalently the subbuilding fixed by $G$ at infinity is a rank one building. Let $\Pi$ be the set of vertices of this rank one building. The elements of $\Pi$ will be denoted by roman letters, elements of $\Lambda$ by greek letters.

Let us recapitulate what we already had obtained in Section 5.2 of Part I (\cite{Muh-Str-Mal:*}). We there constructed a subset $\Lambda'$ of the completion of $\Lambda$, together with a set of injections from an Euclidean space (which is in the present case the real line) to this set. Each apartment of the fixed building at infinity (so for every two different elements $p,q \in \Pi$) there is one such embedding of the real line in $\Lambda'$ (so in fact a geodesic line), which we will denote by $]p,q[$. As $\Lambda'$ satisfies Condition (A4), it follows that if $p,q$ and $r$ are three vertices in $\Pi$, that then the geodesic lines $]p,q[$ and $]p,r[$ share a closed half-line. Let $(qpr)$ denote the endpoint of this half-line. 

Given three pairwise different vertices $p,q$ and $r$ in $\Pi$, the geodesic lines associated with these vertices define a (possibly trivial) triangle with corners $(qpr)$, $(prq)$ and $(rqp)$, which we denote by $\triangle pqr$.

We denote closed line segments in $\Lambda$ with endpoints $\alpha$ and $\beta$ by $[\alpha,\beta]$. Open and half-open segments are written similarly.  We also can consider the case where one point is at infinity. For example $[\alpha,p[$, with $p \in \Pi$, is the center geodesic ray of the sector-face based at $\alpha$ whose direction is the stabilized simplex of $\Lambda_\infty$ corresponding with $p$. In particular when $\alpha$ lies on a geodesic line $]p,q[$ ($p,q \in \Pi$), then $[\alpha,p[$ is a closed half-line of this line.



Let $S_\infty$ be some non-maximal simplex at infinity. We now define a \emph{relative distance $\dd_{S_\infty}$ with respect to $S_\infty$} on pairs of points in $\Lambda$. Consider the $\R$-building $T(S_\infty)$ as constructed in~\cite[\S 4]{Muh-Str-Mal:*}. If $\alpha$ and $\beta$ are two points of $\Lambda$, then the sector-faces $S_\alpha$ and $S_\beta$ correspond to two points of $T(S_\infty)$. The distance between this pair in $T(S_\infty)$ will be our relative distance $\dd_{S_\infty}(\alpha,\beta)$. If $S_\infty$ is a chamber of $\Lambda_\infty$ we define the relative distance $\dd_{S_\infty}$ to be identically zero.

It is clear that the relative distance of a point to itself is zero, that it is symmetric and that the triangle inequality is satisfied. Hence the relative distance forms a pseudo-metric. The relative distance between two points is smaller or equal than the usual distance because of~\cite[Cor. 2.11]{Par:00}.

\subsection{Properties of triangles}\label{section:reldist}

In this section we take a closer look to triangles $\triangle pqr$ with $p,q,r \in \Pi$. We start with a lemma describing the local geometry of the minimal angles.

\begin{lemma}\label{lemma:mangle}
Let $\Delta$ be a weak spherical building of type $(W,S)$ on which a group $H$ acts. Assume that the action of $H$ admits a Tits diagram $(M,\Gamma,A')$ such that $A' \subset A$ (note that the same diagram automorphisms are involved as those arising from the action of $G$ on $\Lambda_\infty$). Let $P$ and $Q$ be two stabilized simplices of $\Delta$ with type the isotropic orbit of the Tits diagram $(M,\Gamma,A)$, such that their centers lie on the minimal angle $\pi/3$ of this Tits diagram. Then the convex hull of both centers in the geometric realization is a line segment of which the two halves can each be covered by one chamber.
\end{lemma}

\proof
One can assume that $\Gamma$ is trivial. If this is not the case one can reduce the problem to the type-preserving case in the same way as discussed in Sections 2.3 and 2.6 of Part I.

As $\Gamma$ is assumed to be trivial, $P$ and $Q$ consist respectively of a single vertex $p$ and $q$. As $p$ and $q$ are not opposite, one has that the convex hull of both is a line segment $L$ of an apartment (so the Coxeter complex geometrically realized on a sphere) $\Sigma$ of $\Delta$.

We say that this segment $L$ \emph{crosses} some wall of $\Sigma$ if it intersects $L$ and does not contain $p$ or $q$. Suppose that this is indeed the case for a wall $M$. The image $q'$ of $q$ under the reflection of $\Sigma$ associated to $M$ is again the center of a simplex of the same type as $P$ and $Q$, but the angle between $p$ and $q'$ is strictly smaller than  the angle between $p$ and $q$. This is only possible if $p=q'$. In particular this implies that $L$ crosses at most one wall. Moreover there is exactly one wall crossed because if there was no wall crossed then $P$ and $Q$ would lie in the same chamber, which is impossible (as $P$ and $Q$ have the same type).

Let $M$ be this unique wall. The two parts of $L$ divided by this wall (notice these two part have equal length as the associated reflection has to map $q$ to $p$) are each covered by a single chamber as no walls are crossed in a single part of $L$.
\qed

We now take a look at some metric properties of the triangles.

\begin{lemma}\label{lemma:tri} Suppose the triangle $\triangle pqr$ is not trivial (i.e., no two corners are the same), then it is equilateral and the corners have angles equal to $\pi/3$. Moreover the convex hull of the triangle is \emph{flat}, i.e. it is isometric to the convex hull of an equilateral triangle in the real Euclidean plane.
\end{lemma}
\proof
By the second assumption of the main result the minimal angle equals $\pi/3$. Using the natural type-preserving epimorphism from the building at infinity to the residues this implies that the angles in the corners have to be at least $\pi/3$. As the sum of the angles in all three corners in a triangle of a CAT(0)-space  has to be less than or equal to $\pi$ (see~\cite[Prop. II.1.7(4)]{Bri-Hae:00}), one knows that the angles are exactly $\pi/3$. Equilaterality and flatness then follow directly from \cite[Prop. II.2.9]{Bri-Hae:00}.
\qed

\begin{lemma}\label{lemma:tri2}
Let $p,q$ and $r$ be three different vertices in $\Pi$. Let $\alpha$ be a point of the geodesic line $]q,r[$, then the relative distance $\dd_p((qpr),\alpha)$ equals
\begin{itemize}
\item
zero if $\alpha \in {} ] r,(prq) ] \cup [(pqr),q[$,
\item
the minimum of $\frac{\sqrt{3}}{2}  \dd((pqr),\alpha)$ and $\frac{\sqrt{3}}{2} \dd((prq),\alpha)) $ if $\alpha \in [(pqr),(prq)]$.
\end{itemize}
\end{lemma}
\proof
We first consider the possibility that $\alpha \in {} ] r,(prq) ] \cup [(pqr),q[$, or more generally that $\alpha$ lies on one of the geodesic  lines $]p,q[$ or $]p,r[$. All of the sector-faces with direction $p$ based on points of  both of these geodesic lines are in the same asymptotic class (as every two such sector-faces share at least a geodesic ray). So $\dd_p((qpr),\alpha)=0$.

We now handle the second case. Note that this can only occur when the triangle $\triangle pqr$ is not trivial. From Lemma~\ref{lemma:tri} it follows that the angle between the half-lines $[(qpr),q[$ and $[(qpr),r[$ in the point $(qpr)$ is $\pi/3$.  By symmetry reasons one can assume that $\alpha \in [(pqr),\beta]$ where $\beta$ is the midpoint of the segment $[(pqr),(prq)]$. Lemma~\ref{lemma:mangle} applied to the residue of the point $(qpr)$ and the simplices in this residue corresponding  to the half-lines $[(qpr),q[$ and $[(qpr),r[$ (note that the projections of $q$ and $r$ on this residue are exactly the centers of these simplices) implies that a neighborhood of $(qpr)$ in the convex hull of the triangle $\triangle pqr$ is covered by two germs of sectors. In fact the result also implies that some neighborhood of $(qpr)$ in the convex hull of the triangle with corners $(qpr)$, $(pqr)$ and $\beta$ is covered by one sector $S$.

By Lemma 2.2 of Part I there exists an apartment $\Sigma$ of $(\Lambda,\cF)$ containing both the germ of $S$ and the point $\alpha$. This apartment contains a sector $S'$ with the same germ as $S$ and containing $\alpha$. As $S'$ has the same germ as $S$, it covers an initial part of the geodesic ray $[(qpr),q[$. Because the geodesic ray $[(qpr),p[$ is opposite to this ray in the point $(qpr)$, one can find a sector $S''$ based at $(qpr)$ such that $S''_\infty$ contains $p$ (so $S''$ contains the geodesic ray $[(qpr),p[$) such that its germ is opposite the germ of $S'$ in the residue at $(qpr)$. Lemma 2.5 of Part I implies that the sectors $S'$ and $S''$ lie in a unique apartment $\Sigma'$ of $(\Lambda,\cF)$. Considering this apartment which contains the points $\alpha$, $(qpr)$ and at infinity $p$, one calculates that $\dd_p((qpr),\alpha)$ equals $\sin(\pi/3) \dd((pqr),\alpha) =  (\sqrt{3}/2) \dd((pqr),\alpha)$.
\qed


\subsection{Equivalences and retractions}
In this section we construct a quotient $\widetilde\Lambda$ of $\Lambda'$. We do this by defining an equivalence relation, which in its turn is defined by combinations of elementary equivalences. The set $\widetilde\Lambda$ will turn out to be the point set of the desired $\R$-tree as in the statement of the main result.

For each three pairwise different vertices $p$, $q$ and $r$ in $\Pi$, we define elementary equivalences of the points of the triangle $\triangle pqr$ as follows: a point $\alpha$ on the segment $[(qpr),(prq)]$ such that $\dd((qpr),\alpha) \leq \dd((prq),\alpha)$ is elementary equivalent with the point $\beta$ on the segment $[(qpr),(pqr)]$ with the same distance to the corner $(qpr)$ as $\alpha$. One can visualize these elementary equivalences in $\triangle pqr$ as follows.

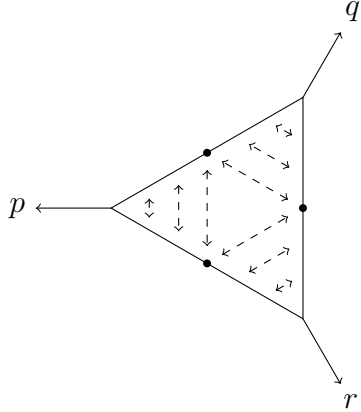
\begin{center}
\begin{tikzpicture}
\draw (60:1.7) -- (180:1.7) -- (-60:1.7) -- cycle;
\draw[->] (60:1.7) --  (60:2.7) node[anchor=240] {$q$};
\draw[->] (-60:1.7) --  (-60:2.7) node[anchor=120] {$r$};
\draw[->] (180:1.7) --  (180:2.7) node[anchor=0] {$p$};

\fill (0:0.85) circle (1.5pt);
\fill (120:0.85) circle (1.5pt);
\fill (240:0.85) circle (1.5pt);

\draw[<->, dashed] (0+10:.66) -- (120-10:.66);
\draw[<->, dashed] (0+39:0.86) -- (120-39:0.86);
\draw[<->, dashed] (54:1.2) -- (120-54:1.2);
\draw[<->, dashed] (120+10:0.66) -- (240-10:0.66);
\draw[<->, dashed] (120+39:0.86) -- (240-39:0.86);
\draw[<->, dashed] (120+54:1.2) -- (240-54:1.2);
\draw[<->, dashed] (240+10:0.66) -- (0-10:0.66);
\draw[<->, dashed] (240+39:0.86) -- (0-39:0.86);
\draw[<->, dashed] (240+54:1.2) -- (0-54:1.2);


%
%
\end{tikzpicture}
\end{center}

Combining these elementary equivalences for all choices of $p$, $q$ and $r$ we obtain an equivalence relation $\sim$ on points of $\Lambda'$ (so two points are equivalent if you can get from one point to the other using a finite number of elementary equivalences). Let $\widetilde\Lambda$ be the quotient of $\Lambda'$ by this equivalence relation.

The composition of the quotient map from $\Lambda'$ to $\widetilde\Lambda$ with the set of charts $\cF'$ to $\Lambda'$ defines a set of charts $\widetilde\cF$ to $\widetilde\Lambda$. In order to prove injectivity of these charts we have to show that for each two points $\alpha$ and $\beta$ on $]p,q[$ with $p,q \in \Pi$ it holds that $\alpha \sim \beta$ if only if $\alpha = \beta$. For this purpose we define a retraction $\rho_{p,q}$ from $\Lambda'$ to $]p,q[$. Let $\Sigma$ be an apartment of $\Lambda$ containing $]p,q[$ and $S$ a sector of $\Sigma$ containing $p$ at infinity, then there exists a retraction $\rho_{S,\Sigma}$ (see~\cite[Prop. 1.20]{Par:00}) mapping $\Lambda$ to $\Sigma$. Let $\alpha$ be a point of $\Lambda'$ and let $\alpha'$ be the orthogonal projection in $\Sigma$ of the image of $\alpha$ under $\rho_{S,\Sigma}$ on $]p,q[$. We now define $\rho_{p,q}(\alpha)$ to be the point on $]p,q[$ at distance $\dd_p(\alpha,\beta) / \sqrt{3}$ from $\alpha'$ towards $p$, where $\beta$ is any point on $]p,q[$ (note that this relative distance is independent of the choice of $\beta$).

The following lemma assures that the retraction $\rho_{p,q}$ is well-defined.

\begin{lemma}
This definition is independent of the choice of $\Sigma$ or $S$.
\end{lemma}
\proof
The choice of $\Sigma$ and $S$ only matter for the projection. Interpreting this projection as a function to $\R$ and using the definition of the retraction $\rho_{S,\Sigma}$ one observes that we are implicitly defining a Busemann function (see~\cite[Def. II.8.17]{Bri-Hae:00}) associated to a geodesic ray with direction $p$. As such Busemann functions only differ by an additive constant (\cite[Ex. II.8.23(1)]{Bri-Hae:00}), and as the projection is constant on $]p,q[$, this construction is independent of the choice of $\Sigma$ or $S$.\qed




Another way to define this retraction is as follows: consider a geodesic line $]r,s[$, let $\beta$ be the midpoint of the segment $[(prs), (psr)]$. This point divides the geodesic line $]r,s[$ in two half-lines. Without loss of generality consider the half-line $[\beta,r[$. First we map this ray to a half-line of the geodesic line $]p,r[$ by `combing' from $r$, and then we `comb' from $p$ to map the half-line to a half-line of $]p,q[$. (By `combing' of one geodesic line to another geodesic line sharing a half-line with the first we mean applying the unique isometric map preserving the intersection.) 

Observe that for both definitions one has that $\rho_{p,s}$ and $\rho_{p,q}$ are identical up to a clear isometry from $]p,s[$ to $]p,q[$. Using this observation and Lemmas~\ref{lemma:tri} and~\ref{lemma:tri2} it follows that these two definitions are equivalent. We will not make use of the first definition later on, its only purpose is to ensure that the second definition is well-defined (which follows from the equivalence).

We will return to the problem of injectivity in Section~\ref{sec:inj}.


\subsection{Possible configurations on four ends}\label{section:fourends}
In this section we investigate which configurations of 6 geodesic lines  with 4 vertices as ends are possible (up to permutation of these ends), which we will use to study the retractions. Let $p$, $q$, $r$ and $s$ be four pairwise different vertices in $\Pi$.  These geodesic lines form 4 (equilateral) triangles, one can suppose without loss of generality that the triangle $\triangle pqr$ has the longest sides. (Note that there may be other triangles with the same longest side lengths.)


The first main case we consider is where the point $(qps)$ lies on $]p,(qpr)[$. This implies that $(qps)=(rps)$. Because the triangle $\triangle pqr$ has the longest sides one has that the point $(pqs)$ lies on $]p,(pqr)[$. One subcase now is that $(pqs)$ lies on $]p,(qpr)[$, then the triangle $\triangle pqr$ equals $\triangle sqr$, and $\triangle psr$ equals $\triangle psq$. So we end up in the following situation:

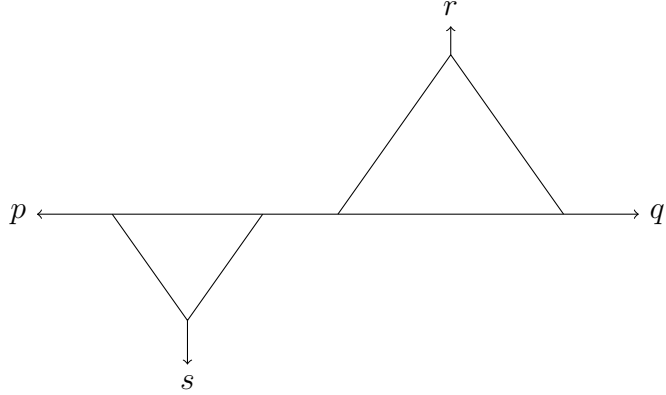
\begin{center}
\begin{tikzpicture}

\draw[<->] (-4,0) node[anchor=east] {$p$} -- (4,0) node[anchor=west] {$q$}; 
\draw (-3,0) -- (-2,-1.414) --  (-1,0);
\draw[->] (-2,-1.414) --  (-2,-2) node[anchor=north] {$s$};

\draw (0,0) -- (1.5,2.121) --  (3,0);
\draw[->] (1.5,2.121) --  (1.5,2.5) node[anchor=south] {$r$};

\end{tikzpicture}
\end{center}

For the second subcase assume that $(pqs)$ lies on $[(qpr),(pqr)[$. A first observation is that $(pqr)=(rqs)$. Also, because the triangle $\triangle pqr$ has longer or equal sides than the triangle $\triangle qrs$ one has that $(qsr) \in {} ](psq),(pqr)[$. This implies that $(psq) =(psr)$, and combined with $(qps) =(rps)$ that the triangles $\triangle pqs$ and $\triangle prs$ have the same side lengths. Hence the point $(prs)$ lies on $[(qpr),(prq)[$. Implications of this are that $(prq)=(qrs)$ and that the triangles $\triangle pqr$ and $\triangle qrs$ have the same side lengths. Finally we remark that because of this $(qsr)$ lies on $](psq),(pqs)]$. The following diagram summarizes the situation (where the triangles are all equilateral, and the points $(prs)$, $(pqs)$, $(qpr)$ and $(qsr)$ lie pairwise on the same distance).

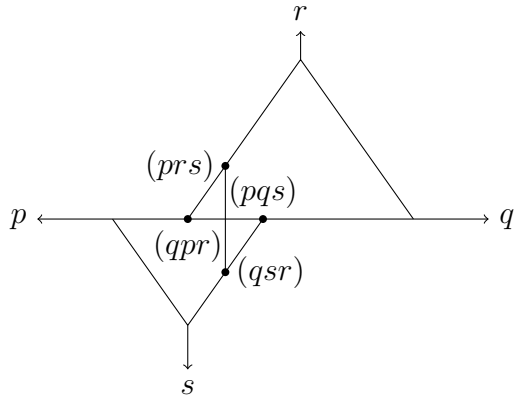
\begin{center}
\begin{tikzpicture}

\draw[<->] (-2,0) node[anchor=east] {$p$} -- (4,0) node[anchor=west] {$q$}; 
\draw (-1,0) -- (-0,-1.414) --  (1,0);
\draw[->] (-0,-1.414) --  (-0,-2) node[anchor=north] {$s$};
\draw (0.5,-.707) -- (0.5,.707);
\fill (0.5,.707) circle (1.5pt) node[anchor=east] {$(prs)$};
\fill (0,0) circle (1.5pt) node[anchor=north] {$(qpr)$};
\fill (1,0) circle (1.5pt) node[anchor=south] {$(pqs)$};

\fill (0.5,-.707) circle (1.5pt) node[anchor=west] {$(qsr)$};
\draw (0,0) -- (1.5,2.121) --  (3,0);
\draw[->] (1.5,2.121) --  (1.5,2.5) node[anchor=south] {$r$};

\end{tikzpicture}
\end{center}


The second main case is when the points $(pqs)$ and $(qps)$ lie on $[(pqr),(qpr)]$ and that $(pqr)= (rqs)$, $(qpr)= (rps)$ (the two last conditions are automatically fulfilled unless we are in the limit case $(pqr) =(pqs)$ or $(qps) = (qpr)$). Note that due to the maximality of the triangle $\triangle pqr$ the points $(prs)$ and $(qrs)$ cannot lie on $](prq),r[$. If the point $(prs)$ would lie on $[(qpr),(prq)[$, then the point $(qrs)$ has to equal $(prq)$. So by interchanging $p$ and $q$ we can assume without loss of generality that $(prs) =(prq)$. It follows that the triangles $\triangle prq$ and $\triangle prs$ have the same side lengths. This implies that $(psr)$ lies on $[(psq),s[$. So we have the following situation.

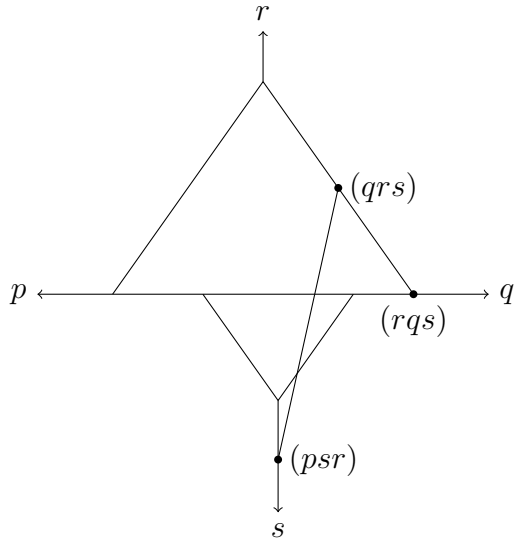
\begin{center}
\begin{tikzpicture}

\draw[<->] (-1,0) node[anchor=east] {$p$} -- (5,0) node[anchor=west] {$q$}; 
\draw (1.2,0) -- (2.2,-1.414) --  (3.2,0);
\draw[->] (2.2,-1.414) --  (2.2,-2.9) node[anchor=north] {$s$};
\fill (4,0) circle (1.5pt) node[anchor=north] {$(rqs)$};
\fill (3,1.414) circle (1.5pt);
\fill (2.2,-2.2) circle (1.5pt);

\draw (3,1.414) node[anchor=west] {$(qrs)$} -- (2.2,-2.2) node[anchor=west] {$(psr)$} ;
\draw (0,0) -- (2,2.828) --  (4,0);
\draw[->] (2,2.828) --  (2,3.5) node[anchor=south] {$r$};

\end{tikzpicture}
\end{center}

Note that it is difficult to picture the situation precisely in two dimensions. Keep in mind that all the triangles are equilateral, for instance the pairwise distances between the points $(qrs)$, $(rqs)$ and $(psr)$ have to be the same. Also remark that the distances $\dd((rqs),(pqs))$ and $\dd((psr),(psq))$ are the same.

The limit configurations mentioned in the second main case reduce to one of the above mentioned cases (by permuting  $p, q, r$ and $s$) except the following configuration: $(rps) \in{} ](qpr),(prq)[$, $(qrs) \in{} ](prq),(pqr)[$ and $(pqs) \in {} ](pqr),(qpr)[$. We will now show that this is impossible. The retraction $\rho_{s,p}$ maps the entire configuration to the geodesic line $]s,p[$. It is easily seen, using the triangle $\triangle prs$, that $(qpr)$ is mapped `(strictly) further from $s$' than $(prq)$, likewise $(prq)$ is mapped further than $(pqr)$, and $(pqr)$ further than $(qpr)$. This would imply that $(qpr)$ is mapped strictly further than $(qpr)$, which is clearly a contradiction.

We conclude that, up to permutation of the ends, only three essentially different configurations (corresponding to the three diagrams above) are possible.

\subsection{Retractions and injectivity}\label{sec:inj}

We now return to problem of showing that the charts $\widetilde\cF$ are injective. Consider the retraction $\rho_{a,b}$. Remark that if two points $\alpha$ and $\beta$ are identified in a triangle of the form $\triangle abc$, then their images under the retraction $\rho_{a,b}$ are the same. We want to prove this for identifications in all other possible triangles. Remark that we only have to consider triangles of the form $\triangle bcd$ (because the retractions of the form $\rho_{a,\cdot}$ are identical up to isometries).

So we have to consider four vertices $a$, $b$, $c$ and $d$. But as all possible configurations of geodesic lines with four vertices as ends are determined in the last section it is straightforward to verify that if two points $\alpha$ and $\beta$ are identified in the triangle $\triangle bcd$ that then their images under $\rho_{a,b}$ are the same. So we conclude that equivalent points have the same image under $\rho_{a,b}$.

In particular this implies that two different points of $]a,b[$ cannot be equivalent as their images are different. This proves injectivity of the charts in $\widetilde\cF$.



\subsection{$(\widetilde\Lambda,\widetilde\cF)$ is an $\R$-tree}
We now verify that $(\widetilde\Lambda,\widetilde\cF)$ is indeed an $\R$-tree, proving the main result.

Conditions (A1) and (A4) are easily seen to be true for $(\widetilde\Lambda,\widetilde\cF)$. Condition (A5) holds as well, this as the identifying of the triangles makes it so that the intersection of three apartments of  $(\widetilde\Lambda,\widetilde\cF)$, pairwise sharing half-lines, is nonempty (as the midpoints of the sides of the triangles are identified with each other).

For Condition (A2) consider two geodesic lines $]a,b[$ and $]c,d[$. Using the possible configurations on the four ends $a,b,c$ and $d$ described in Section~\ref{section:fourends} one checks that a closed segment of each of the two geodesics line is identified, and that the points not on these segments cannot be identified to each other as they can be mapped apart by retractions. So the intersection of two apartments of $(\widetilde\Lambda,\widetilde\cF)$ is closed and convex. Also the metrics agree on both segments, hence (A2) holds.

In order to prove Condition (A3) consider two points $\alpha$ and $\beta$ of $\Lambda$. By construction of $\Lambda$ one can find $a,b,c$ and $d$ in $\Pi$ such that $\alpha \in ]a,b[$ and $\beta \in ]c,d[$. By considering the three cases from Section~\ref{section:fourends} one observes that there exists a pair of points, one equivalent to $\alpha$, the other to $\beta$, both lying on the same geodesic $]a,b[, ]a,c[, ]a,d[, ]b,c[, ]b,d[$ or $]c,d[$. Passing over to the quotient $\widetilde\Lambda$ this proves Condition (A3). Completely analogous one proves the following stronger Condition (A3').

\begin{itemize}
 \item[(A3')] Any two germs are contained in a common apartment.
\end{itemize}

Theorem 1.21 of~\cite{Par:00} states that if Conditions (A1)-(A4) are already satisfied, then Conditions (A3') and (TI) are equivalent, so this proves that $(\widetilde\Lambda,\widetilde\cF)$ is  an $\R$-tree.

This completes the proof of the main result.

\section{Tits diagrams with minimal angle $\pi/3$}\label{subsection:angle}

We now investigate to which diagrams one can apply the main result. Like in Part I, we only need to consider the type-preserving case.

There is a large class of diagrams with minimal angle $\pi/3$. We first give a list of these.

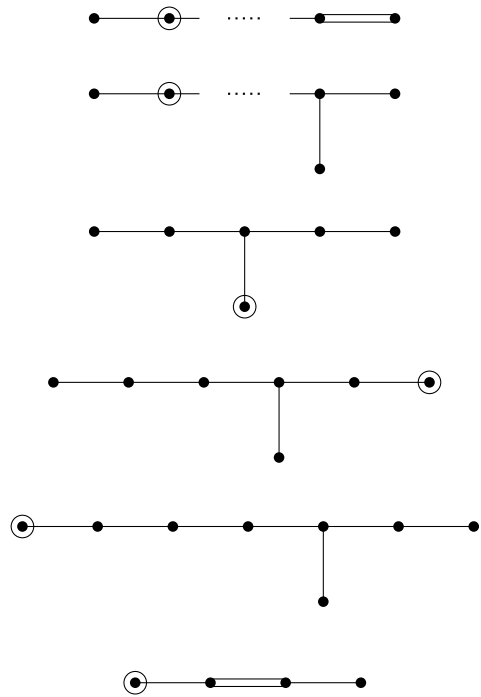
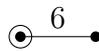
\begin{center}
\begin{tikzpicture}
\fill (-2,0) circle (2pt);
\fill (-1,0) circle (2pt);
\fill (1,0) circle (2pt);
\fill (2,0) circle (2pt);

\draw (-1,0) circle (.15);
\draw[dotted, thick] (-.22,0) -- (.22,0); 

\draw (-2,0) -- (-.6,0); 
\draw (.6,0) -- (1,0); 

\draw (1,.05) -- (2,.05); 
\draw (1,-.05) -- (2,-.05); 

\end{tikzpicture}
\end{center}
\begin{center}
\begin{tikzpicture}
\fill (-2,0) circle (2pt);
\fill (-1,0) circle (2pt);
\fill (1,0) circle (2pt);
\fill (1,-1) circle (2pt);

\fill (2,0) circle (2pt);

\draw (-1,0) circle (.15);
\draw[dotted, thick] (-.22,0) -- (.22,0); 
\draw (-2,0) -- (-.6,0); 
\draw (.6,0) -- (2,0); 

\draw (1,0) -- (1,-1);

\end{tikzpicture}
\end{center}

\begin{center}
\begin{tikzpicture}
\fill (-1,0) circle (2pt);
\fill (0,0) circle (2pt);
\fill (1,0) circle (2pt);
\fill (2,0) circle (2pt);
\fill (3,0) circle (2pt);
\fill (1,-1) circle (2pt);
\draw (1,0) -- (1,-1); 

\draw (-1,0) -- (3,0); 

\draw (1,-1) circle (.15);

\end{tikzpicture}
\end{center}

\begin{center}
\begin{tikzpicture}
\fill (-2,0) circle (2pt);
\fill (-1,0) circle (2pt);
\fill (0,0) circle (2pt);
\fill (1,0) circle (2pt);
\fill (2,0) circle (2pt);
\fill (3,0) circle (2pt);
\fill (1,-1) circle (2pt);
\draw (1,0) -- (1,-1); 

\draw (-2,0) -- (3,0); 

\draw (3,0) circle (.15);

\end{tikzpicture}
\end{center}

\begin{center}
\begin{tikzpicture}
\fill (-3,0) circle (2pt);
\fill (-2,0) circle (2pt);
\fill (-1,0) circle (2pt);
\fill (0,0) circle (2pt);
\fill (1,0) circle (2pt);
\fill (2,0) circle (2pt);
\fill (3,0) circle (2pt);
\fill (1,-1) circle (2pt);
\draw (1,0) -- (1,-1); 

\draw (-3,0) -- (3,0); 

\draw (-3,0) circle (.15);

\end{tikzpicture}
\end{center}

\begin{center}
\begin{tikzpicture}
\fill (-2,0) circle (2pt);
\fill (-1,0) circle (2pt);
\fill (1,0) circle (2pt);
\fill (0,0) circle (2pt);

\draw (-2,0) circle (.15);

\draw (-2,0) -- (-1,0); 
\draw (0,0) -- (1,0); 

\draw (-1,.05) -- (0,.05); 
\draw (-1,-.05) -- (0,-.05); 

\end{tikzpicture}
\end{center}

and
\begin{center}
\begin{tikzpicture}
\fill (-1,0) circle (2pt);
\fill (0,0) circle (2pt);

\draw (-.5,.25) node {6};
\draw (-1,0) -- (0,0); 

\draw (-1,0) circle (.15);

\end{tikzpicture}
\end{center}

The argument for these cases is the following (see also ~\cite[Ch. VI, \S 4.3]{Bou:68}). To these Coxeter systems one can associate a crystallographic root system $\Phi$. The nodes then correspond to a basis of simple roots. If one adds an extra node for the highest root of $\Phi$ then one obtains the \emph{extended Coxeter diagram}. The additional edges then describe the angles between the basis and the highest root: no edge means an angle of $\pi/2$ and a single edge an angle of $\pi/3$. In the cases above the highest root is only connected with a single edge to the encircled node. So the highest root is perpendicular to all roots of the basis but the encircled node. This implies that this highest root is in fact a vertex of the encircled type.

Considering the orbit of this vertex under the Coxeter group one deduces that the vertices of the encircled type correspond exactly to the long roots (`long' in terms of the Euclidean length). As the angle between two roots of the same length in a crystallographic root system can only be $0$, $\pi/3$, $\pi/2$, $2\pi/3$ or $\pi$, one sees that the minimal angle is at least $\pi/3$. From the existence of $\mathsf{A}_2$-subdiagrams in all but the last case it follows that $\pi/3$ indeed occurs as angle. For the last case this is trivial to check.

Alternatively one can use direct calculations to verify these angles.

\begin{remark}\rm
If the Coxeter diagram is not simply laced then there are multiple ways to associate a crystallographic root system to it. However this choice doesn't influence the possible angles so we can make a suitable choice such that the isotropic orbit corresponds to a long root.
\end{remark}

\section{Existence of $\R$-buildings corresponding to exceptional forms}\label{section:existence}
In this section we will discuss the background and proof of the main corollary for algebraic groups and especially exceptional groups. 

We start by some notations and previously known facts. In particular we introduce the framework for Galois descent as developed in~\cite{Bru-Tit:84}. Let $G(K)$ be an absolutely simple algebraic group of relative rank at least one. To such a group there corresponds a spherical building $\Delta_K(G)$. If $\nu$ is a valuation on $K$, there may correspond an $\R$-building $I_{K,\nu}(G)$ to this group (cfr.~\cite[Thm. 4]{Tit:86}). If $\nu$ is a complete valuation then the complete building at infinity of this $\R$-building is the spherical building $\Delta_K(G)$.

The goal is now to prove the conjecture that if $\nu$ is a complete valuation, then such an $\R$-building $I_{K,\nu}(G)$ always exists. So assume that $\nu$ is complete. From~\cite[5.1.4]{Bru-Tit:84} it follows that there exists a finite Galois extension $L$ of $K$ such that the algebraic group $G(L)$ is quasi-split. The valuation $\nu$ extends uniquely to a complete valuation $\omega$ of $L$ by~\cite[Thm. 7.1.1]{Cas:86}.  One can now construct an $\R$-building $I_{L,\omega}(G)$ corresponding with this group and the valuation $\omega$ (see~\cite[\S 4]{Bru-Tit:84}). The Galois group $\mbox{Gal}(L/K)$ acts boundedly on $I_{L,\omega}(G)$, again by~\cite[5.1.4]{Bru-Tit:84}. The fixed structure of $\mbox{Gal}(L/K)$ on the complete building at infinity is the spherical building $\Delta_K(G)$, moreover this action admits a Tits diagram for which the possibilities are listed in~\cite{Tit:66}.

So we are in a situation where we can try to apply the main results from both parts of this two-part paper. Examining the list in~\cite{Tit:66} of exceptional Tits diagrams, there are only four Tits diagrams for which none of the main results apply (so these have minimal angles less than $\pi/3$, or have minimal angle equal to $\pi/3$ but are of relative rank strictly greater than one). These Tits diagrams are:


%
%

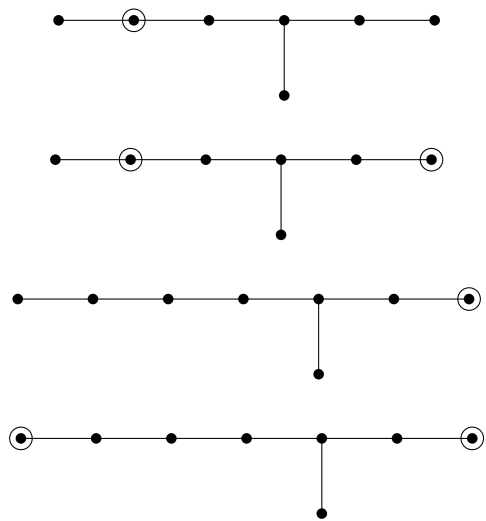
\begin{center}
\begin{tikzpicture}
\fill (-2,0) circle (2pt);
\fill (-1,0) circle (2pt);
\fill (0,0) circle (2pt);
\fill (1,0) circle (2pt);
\fill (2,0) circle (2pt);
\fill (3,0) circle (2pt);
\fill (1,-1) circle (2pt);
\draw (1,0) -- (1,-1); 

\draw (-2,0) -- (3,0); 

\draw (-1,0) circle (.15);

\end{tikzpicture}
\end{center}

\begin{center}
\begin{tikzpicture}
\fill (-2,0) circle (2pt);
\fill (-1,0) circle (2pt);
\fill (0,0) circle (2pt);
\fill (1,0) circle (2pt);
\fill (2,0) circle (2pt);
\fill (3,0) circle (2pt);
\fill (1,-1) circle (2pt);
\draw (1,0) -- (1,-1); 

\draw (-2,0) -- (3,0); 
\draw (-1,0) circle (.15);
\draw (3,0) circle (.15);

\end{tikzpicture}
\end{center}

\begin{center}
\begin{tikzpicture}
\fill (-3,0) circle (2pt);
\fill (-2,0) circle (2pt);
\fill (-1,0) circle (2pt);
\fill (0,0) circle (2pt);
\fill (1,0) circle (2pt);
\fill (2,0) circle (2pt);
\fill (3,0) circle (2pt);
\fill (1,-1) circle (2pt);
\draw (1,0) -- (1,-1); 

\draw (-3,0) -- (3,0); 

\draw (3,0) circle (.15);

\end{tikzpicture}
\end{center}

\begin{center}
\begin{tikzpicture}
\fill (-3,0) circle (2pt);
\fill (-2,0) circle (2pt);
\fill (-1,0) circle (2pt);
\fill (0,0) circle (2pt);
\fill (1,0) circle (2pt);
\fill (2,0) circle (2pt);
\fill (3,0) circle (2pt);
\fill (1,-1) circle (2pt);
\draw (1,0) -- (1,-1); 

\draw (-3,0) -- (3,0); 
\draw (3,0) circle (.15);
\draw (-3,0) circle (.15);

\end{tikzpicture}
\end{center}

In particular this implies the existence of an affine building $I_{K,\nu}(G)$ in all other cases. In the next two sections we deal with the remaining cases.

%

\subsection{The remaining relative rank two cases}

The remaining relative rank two cases correspond to the Moufang quadrangles of exceptional type $\mathsf{E}_7$ and  $\mathsf{E}_8$. So let $\Gamma$ be a Moufang quadrangle of exceptional type $\mathsf{E}_i$ ($i=7$ or $8$), defined over a field $K$ with a complete valuation $\nu: K \to \mathbb{R} \cup \{\infty\}$, and let $G(K)$ be the corresponding exceptional algebraic group. 

The group $G$ splits over some extension $E$ of $K$. While not important for our proof one can in fact choose $E$ to be a quadratic extension of $K$ (see~\cite{Ste:*}).  Let $\omega$ be the extension of the valuation $\nu$ to $E$. So we can consider the $\R$-building $I_{E,\omega}(G)$ and the Galois involution acting on it. 

We will now prove the existence of an $\R$-building $I_{K,\nu}(G)$ in two ways: one using the underlying algebraic parameter systems and one using the embedding result of Guy Rousseau.

\subsubsection{An algebraic proof}

Applying the methods from Section 5.2 of~\cite{Muh-Str-Mal:*} to $I_{E,\omega}(G)$ and the Galois involution acting on it, one obtains a subset $\Lambda'$ of the points of $I_{E,\omega}(G)$. We make the special choice of fixed point as outlined in Remark 5.5 of~\cite{Muh-Str-Mal:*}. This has as advantage that isometries of $I_{E,\omega}(G)$ centralizing the Galois involution stabilize $\Lambda'$.

Although our methods did not prove that this was an $\R$-building (as we are in the `minimal angle equal to $\pi/3$' case), it is however possible to consider flats in $\Lambda'$ for which the directions are the apartments of the fixed structure $\Delta_K(G)$ at infinity. Pick one such flat $\Sigma'$ of $\Lambda'$. For the corresponding apartment of $\Delta_K(G)$ we consider four consecutive root groups $U_j$ (with $j \in \{1,2,3,4\}$) associated to it, so that they are of the form as described in~\cite[Ex. 16.6]{Tit-Wei:02}.

These root groups stabilize $\Lambda'$ so they map $\Sigma'$ to another flat sharing a half-space with it.   These root groups live in rank one residues of $\Gamma$, for which the corresponding Tits subdiagrams of relative rank one has minimal angles at least $\pi/3$. So we can apply the main results of this two-part paper and prove existence for these $\R$-trees. From these one constructs valuation-like maps $\phi_j: U_j \to \mathbb{R} \cup \{ \infty \}$ (with $j \in \{1,2,3,4\}$) such that each superlevel set with respect to these functions forms a subgroup of the root group, this as in~\cite[Def. 13.8]{Wei:09}. 

By the geometric nature of these methods these maps can be interpreted geometrically. An important consequence is that this allows us to consider the maps $m_\Sigma(u)$ (see~\cite[Def. 3.8]{Wei:09}) and their interactions with the maps $\phi_j$ ($j \in \{1,2,3,4\}$) as in~\cite[Prop. 3.31]{Wei:09}. This allows one to derive the explicit form (up to equipollence and rescaling) of the map $\phi_1$ with the same methods as the analogous result~\cite[Thm. 21.27]{Wei:09}:
\begin{equation*}
\phi_1(x_1(a,t)) = \nu(q(\pi(a) +t)) .
\end{equation*}
If the valuation $\nu$ is discrete and normalized such that $\nu(K^*) = \mathbb{Z}$, then the image of this function is either $\mathbb{Z}$ or $2\mathbb{Z}$ (\cite[Prop. 21.24]{Wei:09}). As the existence problem for the subdiagrams of relative rank one is solved, we have that (because the superlevel sets w.r.t. $\phi_1$ forms a subgroup)
\begin{equation*}
\forall (a,t), (b,s) \in S:  \phi_1(x_1(a,t)), \phi_1(x_1(b,s)) \geq 0 \Rightarrow \phi_1(x_1(a,t)\cdot x_1(b,s) ) \geq 0,
\end{equation*}
or using the explicit form of $\phi_1$ (where $S$ is used to index $U_1$ as in~\cite[(16.6)]{Tit-Wei:02})
\begin{equation*}
\forall (a,t), (b,s) \in S:  \nu(q(\pi(a) + t)), \nu(q(\pi(b) + s)) \geq 0 \Rightarrow \nu(q(\pi(a+b) + t + s + g(b,a))) \geq 0.
\end{equation*}
This can be simplified to the following using~\cite[1.17(C3)]{Wei:06}:
\begin{equation*}
\forall (a,t), (b,s) \in S:  \nu(q(\pi(a) + t)), \nu(q(\pi(b) + s)) \geq 0 \Rightarrow \nu(q(\pi(a)+\pi(b) + t + s + h(b,a))) \geq 0.
\end{equation*}
Note that by completeness of $\nu$ and~\cite[Prop. 19.4]{Wei:09} one has that if $\nu(q(c)) < \nu(q(d))$, that then $\nu(q(c + d)) = \nu(q(c))$. In particular this implies 
\begin{equation*}
\forall (a,t), (b,s) \in S:  \nu(q(\pi(a) + t)), \nu(q(\pi(b) + s)) \geq 0 \Rightarrow \nu(q(h(b,a))) \geq 0.
\end{equation*}
If $k \in K$, then one has that $\nu(q(\pi(ka) + k^2t)) =  4 \nu(k) + \nu(q(\pi(a) + t))$ (see for instance~\cite[Prop. 21.24]{Wei:09}) and that the function $h$ is bilinear over $K$ (see~\cite[Def. 13.18-19]{Tit-Wei:02}). Hence by replacing $(a,t)$ by $(ka,k^2t)$ and $(b,s)$ by $(lb,l^2s)$ with $k,l \in K$ one obtains
\begin{equation*}
\forall (a,t), (b,s) \in S:  \nu(q(\pi(a) + t)) \geq 4m, \nu(q(\pi(b) + s)) \geq 4n \Rightarrow \nu(q(h(b,a))) \geq 2m + 2n,
\end{equation*}
where $m = - \nu(k)$ and $n= -\nu(l)$. If the image of  the valuation $\nu$ is dense in $\R$, then this directly implies Equation~\ref{eq:2} below. We now show this implication in the discrete case. When the valuation $\nu$ is surjective on $\mathbb{Z}$, the above equation is true for all $m,n \in \mathbb{Z}$. As the image of $\phi_1$ is either $\mathbb{Z}$ or $2\mathbb{Z}$, this yields the following inequality:
\begin{equation}\label{ineq}
\forall (a,t), (b,s) \in S:  [\nu(q(\pi(a) + t)) + \nu(q(\pi(b) + s))]/2 - 3 \leq \nu(q(h(b,a))).
\end{equation}

Let $\varpi$ be an uniformizer of the field $K$ (i.e. $\nu(\varpi)=1$), and $p$ an odd prime. Extend the field $K$ to a field $K'$ by adding an element $\beta$ such that $\beta^p = \varpi$. The valuation $\nu$ extends uniquely to a complete valuation $\nu'$ of $K'$ (see~\cite[Thm. 7.1.1]{Cas:86}). Note that the value group of this new valuation contains $\mathbb{Z} / p$. As the degree of the extension is at most $p$, this implies it is exactly $p$, and the value group is $\mathbb{Z} / p$ (see~\cite[\S 1, Lem. 18]{Sch:50}). 

By Springer's theorem (see~\cite[Cor. 18.5]{Elm-Kar-Mer:08}), the anisotropic quadratic form with a norm splitting over $K$  defining $\Gamma$ (for details and definitions see~\cite[\S 12]{Tit-Wei:02}) stays anisotropic when tensored up over $K'$ (and also stays a norm splitting with the same constants). This can be used to embed the exceptional quadrangle $\Gamma$ into a larger exceptional quadrangle $\Gamma'$ defined over $K'$ (and extending $S$ to an $S'$). With the same reasoning as before, one now obtains the stronger version of Equation~\ref{ineq} (because the value group is finer):
\begin{equation*}
\forall (a,t), (b,s) \in S':  [\nu'(q(\pi(a) + t)) + \nu'(q(\pi(b) + s))]/2 - 3/p \leq \nu'(q(h(b,a))).
\end{equation*}
As this is true as well if one restricts to $S$, and for arbitrary large odd primes $p$, one obtains:
\begin{equation}\label{eq:2}
\forall (a,t), (b,s) \in S:  [\nu(q(\pi(a) + t)) + \nu(q(\pi(b) + s))]/2 \leq \nu(q(h(b,a))).
\end{equation}
This implies the existence of an $\R$-building $I_{K,\nu}(G)$ for the exceptional Moufang quadrangle $\Gamma$ defined over $K$ by~\cite[Thm. 21.27]{Wei:09}.

This settles the existence for the two remaining rank two cases.

\subsubsection{Alternative proof using descent}
In this section we show how to prove the existence for the two remaining rank two cases using descent. In particular this approach is based on re-evaluating the choice of fixed point in Section 5.1 of Part I using an embedding result of Guy Rousseau (see Theorem 2.11 of Part I). We discuss this for the Moufang quadrangle of exceptional type $\mathsf{E_7}$. Once this case is settled, the $\mathsf{E}_8$ case is proven completely analogously.

Let $S_\infty$ be a maximal stabilized simplex at infinity. Let $S'_\infty \subset S_\infty$ be a subsimplex of rank one such that deleting the associated isotropic orbit from the Tits diagram yields the following Tits diagram with minimal angle $\pi/3$. (For the $\mathsf{E}_8$ case one needs to delete the isotropic orbit such that one gets an $\mathsf{E}_7$ diagram.)

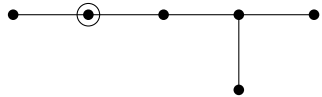
\begin{center}
\begin{tikzpicture}
\fill (-2,0) circle (2pt);
\fill (-1,0) circle (2pt);
\fill (0,0) circle (2pt);
\fill (1,0) circle (2pt);
\fill (2,0) circle (2pt);
\fill (1,-1) circle (2pt);
\draw (1,0) -- (1,-1); 

\draw (-2,0) -- (2,0); 
\draw (-1,0) circle (.15);

\end{tikzpicture}
\end{center}

Let us recapitulate the proof of the existence problem for this Tits diagram. In Section 5.1 of Part I we had to pick a fixed point of a group $H$ in a convex subset $F$ of $T(S_\infty)$ (this set is also bounded by Remark 5.5 of Part I). The choice of this fixed point yielded a `tree with triangles' from which we then constructed a tree. 

Note that classical results (see~\cite{Bru-Tit:72}) handle each Tits diagram of type $\mathsf{D}_6$ possible. (For the $\mathsf{E}_8$ case the combination of both main results of this two-part paper and the previous result for the $\mathsf{E}_7$ case handle each Tits diagram of type $\mathsf{E}_7$ listed in~\cite{Tit:66}.) This implies that we can apply Theorem 2.11 of Part I, and embed this tree nicely in (the completion of) $T(S'_\infty$). This embedding corresponds with a certain choice of fixed point in $F$. The set of all points in $T(S_\infty)$ giving rise to such an embedding of the tree is a bounded and convex subset $F'$ of $F$ (see Lemma 5.8 of Part I).

We now return to the existence problem for the Tits diagram  

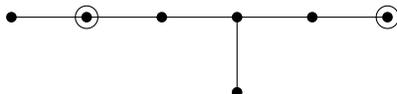
\begin{center}
\begin{tikzpicture}
\fill (-2,0) circle (2pt);
\fill (-1,0) circle (2pt);
\fill (0,0) circle (2pt);
\fill (1,0) circle (2pt);
\fill (2,0) circle (2pt);
\fill (3,0) circle (2pt);
\fill (1,-1) circle (2pt);
\draw (1,0) -- (1,-1); 

\draw (-2,0) -- (3,0); 
\draw (-1,0) circle (.15);
\draw (3,0) circle (.15);

\end{tikzpicture}
\end{center}

Note that we are again dealing with the choice of a fixed point in the same set $F$ in $T(S_\infty)$, this time for some group $H'$. This group $H'$ is contained in the stabilizer of $S_\infty$ in the little projective group of the building $\Delta_K(G)$. This stabilizer of $S_\infty$ also stabilizes the bounded and convex subset $F'$. So instead of picking a fixed point in the set $F$ we may use the center of $F'$ as choice of fixed point.

Following the proof of the main result of Part I this choice yields a subset $\Lambda'$ of the point set of the $\R$-building $I_{E,\omega}(G)$ which is stabilized by the group $G(K)$. The only obstacle for this proof to fully work are possible triangle configurations as discussed in Section 5.2.1 of Part I. As mentioned in that section, we can reduce this problem to the subdiagrams of relative rank one. The problematic subdiagram is the subdiagram of type $\mathsf{D}_6$ pictured above. However the special choice of fixed point in $F'$ guarantees that there are no triangle configurations possible where one of the flats contains $S_\infty$ at infinity. As $G(K)$ acts transitive on all flats of which $\Lambda'$ consists, there are no triangle configurations possible at all. So the proof of the main result of Part I (with the above adaptations) applies to this case.

This gives an alternative answer to the existence problem for the Moufang quadrangle of exceptional type $\mathsf{E}_7$. The $\mathsf{E}_8$ case can now be proven the same way.

\subsection{The remaining relative rank one cases}

Both of the remaining relative rank one cases have an anisotropic kernel containing $\mathsf{D}_i$ with $i$ odd as a connected component. Let us make the following observation.

\begin{lemma}
This anisotropic kernel can be split by a tower of quadratic extensions. 
\end{lemma}
\proof
An algebraic group $H(K)$ of (absolute) type $\mathsf{D}_i$ over a field $K$ acts on a non-singular polar space (which is in fact the associated spherical building $\Delta_K H$) defined by some quadratic form in the projective space $\PG(K,2i-1)$. For a field extension $K'$ of $K$, the group $H(K')$ is split if and only if the building $\Delta_{K'} H$ has rank $i$. Interpreting $\Delta_{K'} H$ as a polar space, this is equivalent to saying that this non-singular polar space contains subspaces of (the maximal possible) dimension $i-1$.

We prove this by induction on $i$. We start with $i=1$, so we are dealing with a non-singular polar space in $\PG(K,1)$. The quadratic form defining the polar space here is a homogeneous quadratic equation in two variables, which has a solution for a suitable quadratic extension $K'$ of $K$. This solution corresponds to a point of the polar space. Note that this argument also holds for singular polar spaces.

Now consider a general $i>1$ and suppose the statement is proven for smaller values. If we restrict to a line of this polar space we can apply the argument for $i=1$ and conclude that there is a quadratic field extension for which the polar space has at least a point. So without loss of generality we can assume that the polar space contains a point $x$. The intersection of this polar space and the tangent hyperplane at  $x$ of this polar space yields a cone with apex $x$ and base a non-singular polar space in $\PG(K,2i-3)$. The induction hypothesis implies that this last polar space, considered over a field $K'$ obtained by a tower of quadratic extensions starting from $K$, contains subspaces of dimension $i-2$. The span of such a subspace and the point $x$ now forms a desired subspace of dimension $i-1$. \qed


Suppose that this anisotropic part would split completely at one such extension, then this would geometrically amount to an involution on a spherical building of type $\mathsf{D}_i$ without fixed points (so every simplex has to be mapped to an opposite one). This involution would be type-preserving as it  arises from an involution of an $\mathsf{E}_7$ or $\mathsf{E_8}$ building. However such an involution cannot be fixed point free as the opposition involution for $\mathsf{D}_i$ is only type-preserving if $i$ is even. We conclude that we cannot split the form completely with one step.

The Tits diagram of an intermediate form has to be one of the following by Tits' list~\cite{Tit:66}.

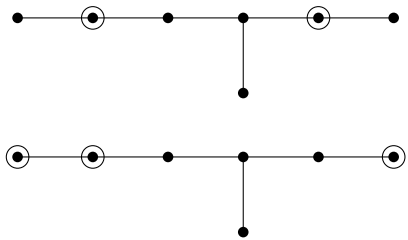
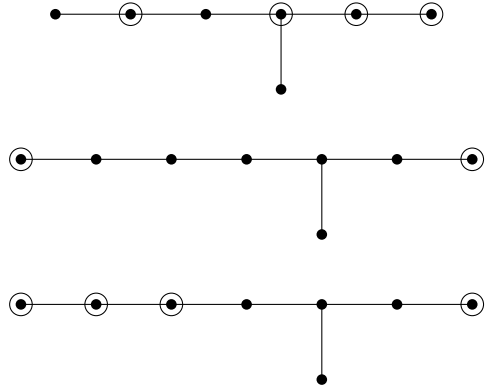
\begin{center}
\begin{tikzpicture}
\fill (-2,0) circle (2pt);
\fill (-1,0) circle (2pt);
\fill (0,0) circle (2pt);
\fill (1,0) circle (2pt);
\fill (2,0) circle (2pt);
\fill (3,0) circle (2pt);
\fill (1,-1) circle (2pt);
\draw (1,0) -- (1,-1); 

\draw (-2,0) -- (3,0); 
\draw (-1,0) circle (.15);

\draw (2,0) circle (.15);

\end{tikzpicture}
\end{center}

\begin{center}
\begin{tikzpicture}
\fill (-2,0) circle (2pt);
\fill (-1,0) circle (2pt);
\fill (0,0) circle (2pt);
\fill (1,0) circle (2pt);
\fill (2,0) circle (2pt);
\fill (3,0) circle (2pt);
\fill (1,-1) circle (2pt);
\draw (1,0) -- (1,-1); 

\draw (-2,0) -- (3,0); 
\draw (-1,0) circle (.15);
\draw (-2,0) circle (.15);

\draw (3,0) circle (.15);

\end{tikzpicture}
\end{center}
\begin{center}
\begin{tikzpicture}
\fill (-2,0) circle (2pt);
\fill (-1,0) circle (2pt);
\fill (0,0) circle (2pt);
\fill (1,0) circle (2pt);
\fill (2,0) circle (2pt);
\fill (3,0) circle (2pt);
\fill (1,-1) circle (2pt);
\draw (1,0) -- (1,-1); 

\draw (-2,0) -- (3,0); 
\draw (-1,0) circle (.15);
\draw (1,0) circle (.15);
\draw (2,0) circle (.15);
\draw (3,0) circle (.15);

\end{tikzpicture}
\end{center}

\begin{center}
\begin{tikzpicture}
\fill (-3,0) circle (2pt);
\fill (-2,0) circle (2pt);
\fill (-1,0) circle (2pt);
\fill (0,0) circle (2pt);
\fill (1,0) circle (2pt);
\fill (2,0) circle (2pt);
\fill (3,0) circle (2pt);
\fill (1,-1) circle (2pt);
\draw (1,0) -- (1,-1); 

\draw (-3,0) -- (3,0); 

\draw (-3,0) circle (.15);
\draw (3,0) circle (.15);

\end{tikzpicture}
\end{center}

\begin{center}
\begin{tikzpicture}
\fill (-3,0) circle (2pt);
\fill (-2,0) circle (2pt);
\fill (-1,0) circle (2pt);
\fill (0,0) circle (2pt);
\fill (1,0) circle (2pt);
\fill (2,0) circle (2pt);
\fill (3,0) circle (2pt);
\fill (1,-1) circle (2pt);
\draw (1,0) -- (1,-1); 

\draw (-3,0) -- (3,0); 
\draw (3,0) circle (.15);
\draw (-3,0) circle (.15);
\draw (-2,0) circle (.15);
\draw (-1,0) circle (.15);

\end{tikzpicture}
\end{center}

In particular this means that one can use Tits diagrams of type $\mathsf{B}_2$, $\mathsf{C_3}$ or $\mathsf{F}_4$ to study the two problematic rank one forms. This intermediate step allow us to apply our geometric methods twice to cases we can already solve. This answers the existence problem for the two remaining rank one cases.


\begin{thebibliography}{99}
\bibitem{Bou:68}
N. Bourbaki, {\em Groupes et alg\`ebres de Lie. Chap. 4--6}, 1968.

\bibitem{Bri-Hae:00}
M. Bridson and A. Haefliger, {\em Metric spaces of non-positive curvature}, Grundlehren Math. Wiss., {\bf 319}, Springer-Verlag, Berlin, 1999.


\bibitem{Bru-Tit:72} 
F. Bruhat and J. Tits, Groupes r\'eductifs sur un corps local: 
I.~Donn\'ees radicielles valu\'ees, {\em Inst.\ Hautes
\'Etudes Sci.\ Publ.\ Math}.\ {\bf 41} (1972), 5--252.

\bibitem{Bru-Tit:84} 
F. Bruhat and J. Tits, Groupes r\'eductifs sur un corps local: 
II.~Sch\'emas en groupes. Existence d'une donn\'ee racidielle valu\'ee, {\em Inst.\ Hautes
\'Etudes Sci.\ Publ.\ Math}.\ {\bf 60} (1984), 5--184.
%

\bibitem{Cas:86}
J. W. S. Cassels, \emph{Local fields}, LMS Student Series \textbf{3}, Cambridge University Press, 1986.

\bibitem{Elm-Kar-Mer:08}
R. Elman, N. Karpenko and A. Merkurjev, \emph{The algebraic and geometric theory of quadratic forms}, Am. Math. Soc., Providence, 2008.

\bibitem{Kna:88}
N. Knarr, Projectivities of generalized polygons, \emph{Ars Combin.} \textbf{25B} (1988), 265--275.

%
%
%

\bibitem{Muh-Str-Mal:*}
B. M\"uhlherr, K. Struyve and H. Van Maldeghem, Descent of affine buildings - I. Large minimal angles, \emph{preprint}.
 
\bibitem{Par:00}
A. Parreau, Immeubles affines: construction par les normes et \'etude des isom\'etries, {\bf in} \emph{Crystallographic Groups and their Generalizations (Kortrijk, 1999)}, Contemp. Math. \textbf{262}, Amer. Math. Soc., Providence, RI,
2000, pp. 263--302.

\bibitem{Rou:77}
G. Rousseau, {\em Immeubles des groupes r\'educitifs sur les corps locaux.}, Th\`ese de doctorat, Publications Math\'ematiques d'Orsay, No. 221-77.68, 1977.


\bibitem{Sch:50}
O. F. G. Schilling, \emph{The theory of valuations}, Mathematical Surveys, no. IV, 1950.

\bibitem{Ste:*}
A. Steinbach, Realizing Moufang quadrangles of type $\mathsf{E}_n$ inside Chevalley groups, \emph{preprint}.

%
%

\bibitem{Tit:66}
J. Tits, Classification of algebraic semisimple groups, {\bf in} \emph{Algebraic groups and Discontinuous Subgroups (Proc. Sympos. Pure Math., Boulder, Colo.,1965)} pp. 33--62 \emph{Amer. Math. Soc., Providence, R.I, 1966}

\bibitem{Tit:86}
J. Tits, Immeubles de type affine, {\bf in} \emph{Buildings and the Geometry of Diagrams, Springer Lecture Notes}
\textbf{1181} (Rosati ed.), Springer Verlag, 1986, pp. 159--190. 

\bibitem{Tit-Wei:02} J.~Tits and R.~Weiss, \emph{Moufang polygons}, Springer-Verlag, 2002.


\bibitem{Wei:06}
R. Weiss, \emph{Quadrangular algebras}, Math. Notes \textbf{46}, Princeton University Press, 2006. 

\bibitem{Wei:09} 
R. Weiss, \emph{The structure of affine buildings}, Annals of Math Studies \textbf{168}, Princeton University Press, Princeton, 2009.

\bibitem{Wei:*}
R. Weiss, On the existence of certain affine buildings of type $\mathsf{E}_6$ and $\mathsf{E}_7$, {\em J. Reine Angew. Math.} \textbf{653} (2011), 135--147.

\end{thebibliography}
\end{document}